\documentstyle[aps,pra,multirow,hhline,epsfig]{revtex}      

\title{Conjugated filter approach for solving Burgers' 
equation with high Reynolds number}

\author{G. W. Wei and Yun Gu \\
 Department of Computational Science,\\ 
National University of Singapore, 
Singapore 117543}

\date{\today}

\begin{document}
\maketitle

\begin{abstract}
 
We propose a conjugated filter 
oscillation reduction scheme for solving  
Burgers' equation with high Reynolds numbers.
Computational accuracy is tested at a moderately 
high  Reynolds 
number for which analytical solution is available.
Numerical results at extremely high Reynolds numbers 
indicate that the proposed scheme 
is efficient, robust and reliable for shock capturing. 
\end{abstract} 
PACS numbers: 02.70.-c, 05.10.-a, 47.40.Nm,52.35.Tc

The fundamental equation for the 
description of complex fluid flow 
is the Navier-Stokes equation, for which
the full solution is still extremely difficult in 
the full domain of physical interest.
Burgers' equation\cite{Burger} is an important simple model 
for the understanding of physical flows.
Simulation of Burgers' equation 
is a natural first step towards developing methods
for computations of complex flows. 
It appears customary to test 
new approaches in computational fluid dynamics 
by applying them to Burgers' 
equation. 
Jamet and Bonnerot solved Burgers' equation by using
isoparametric rectangular space-time finite elements\cite{JamBon}.
Jain and Holla\cite{JaiHol} developed a 
cubic spline approach for coupled Burgers' 
equation.
Varoglu and Finn\cite{VarFin} proposed an isoparametric space-time
 finite element method for solving  Burgers' equation,
utilizing the hyperbolic differential equation associated 
with Burgers' equation. 
They obtain very high accuracy 
and numerical stability with a reasonable number of elements 
and time steps. 
Their method was compared with a least 
square weak formulation of the finite element method by 
Nguyen and Reynen\cite{NguRey}. Caldwell et 
al\cite{Caldwell} further developed the finite element method 
to allow different sizes of the elements at each stage 
based on the feed  back from the previous step. A generalized 
boundary element approach was proposed by Kakuda and 
Tosaka\cite{KakTos}. These authors tabulated their accurate
 results for moderate Reynolds numbers and compared their 
results with those of Varoglu and Finn\cite{VarFin} and 
of Nguyen and Reynen\cite{NguRey}. 
A bidimsional Tau-element method was developed by Ortiz and 
Pun for solving Burgers' equation with accurate results.
Bar-Yoseph et al 
discussed a number of space-time spectral 
element methods for solving  Burgers' equation\cite{BMZY}.
Arina and Canuto\cite{AriCan} treated  Burgers' equation 
by a self-adaptive, domain decomposition method called the 
$\chi$-formulation. Various finite difference schemes for  
Burgers' equation were compared by Biringen and
 Saati\cite{BirSaa}. Recently, Wei et al\cite{weicpc} 
have developed an accurate solver for  Burgers' equation
in one and two space dimensions.
Most recently, Hon and Mao\cite{HonMao} have 
compared performance of  their adaptive 
multiquadric scheme with many other computational 
methods. It is not
our purpose to exhaust the literature.
Despite of much effort, numerical 
solution  of  Burgers' equation is still a nontrivial 
task especially at very high Reynolds numbers
where the 
nonlinear advection leads to shock waves.
In fact, Burgers' inviscid shocks plague many standard
computational algorithms.

The purpose of this communication is to report a novel 
scheme for solving Burgers' equation
for all possible values of Reynolds numbers.
We propose a set of conjugated filters to 
solve the Burgers' equation.
As the first and second derivatives are 
approximated by using two high-pass filters, 
the numerical errors of the  high-pass 
filters at the high frequency region lead to oscillation 
near the Burgers' shock wave front. 
The proposed idea is to effectively eliminate such 
an oscillation by using a conjugated low-pass filter.
This set of high-pass and low-pass filters are 
conjugated in the sense that they are derived from one 
generating function and consequently have essentially the 
same degree of regularity, smoothness, time-frequency 
localization, effective support and bandwidth.   
In the present work, all conjugated filters are 
constructed by using a discrete singular convolution (DSC) 
algorithm\cite{weijcp99}, which is a potential approach for 
dynamical simulation\cite{weiphysica,weijpa20} and numerical 
computations of Hilbert transform and Radon transform.
The DSC low-pass filter is given by 
\begin{equation}\label{K1}
\phi_{\Delta,\sigma;k}(x)=
{\sin{\pi\over\Delta}(x-x_k)\over
{\pi\over\Delta}(x-x_k)}
\exp\left[ {-{(x-x_k)^2\over2\sigma^2}}\right],
\end{equation}     
and its conjugated $n$th order high-pass filters are obtained by  
differentiation   
\begin{equation}\label{deriv}   
\phi^{(n)}_{\Delta,\sigma;k}(x)=
\left({d\over dx}\right)^n 
{\sin{\pi\over\Delta}(x-x_k)\over
{\pi\over\Delta}(x-x_k)}
\exp\left[ {-{(x-x_k)^2\over2\sigma^2}}\right], ~~~n=1,2,...,
\end{equation}
where $\Delta$ is the grid spacing and $\sigma$ is 
a regularization parameter.
FIG. 1 shows the frequency responses of the conjugated DSC low-pass
filter,  1st and 2nd order high-pass filters at $\sigma=3.2\Delta$. 
It is noted that 
all conjugated filters have essentially the same effective bandwidth.
Wavelet multiscale analysis is utilized 
for adaptive oscillation control. 
A detailed theoretical analysis and justification 
of this scheme is accounted elsewhere\cite{GuWei}.

Burgers' equation is  given by 
\begin{equation}\label{Burgers}
{\partial u\over \partial t}+u{\partial u\over\partial x}={1\over {\rm Re}}  
{\partial^2 u\over\partial x^2}, 
\end{equation}                 
where $u(x,t)$ is the dependent variable 
resembling the flow velocity and Re is the 
Reynolds number characterizing the size of viscosity.
The competition between the nonlinear advection and the
 viscous diffusion is 
controlled by the value of Re in  Burgers' equation, and
thus determines the behavior of the solution. 
 We consider Eq. (\ref{Burgers}) using the 
following initial-boundary conditions
\begin{eqnarray}\label{bur2}
u(x,0)&=&\sin(\pi x),\nonumber\\
u(0,t)&=&u(1,t)=0.
\end{eqnarray}
Cole has provided exact solution\cite{Cole} for this problem 
in terms of a series expansion which is readily 
computable roughly for the parameter ${\rm Re}\leq 100$.     
For the parameter ${\rm Re}=100$, the present calculations use 
41 grid points in the interval [0,1]. 
The fourth order Runge-Kutta is used for the time 
integration with a time increment of 0.01.

In the present DSC treatment, 
we use the DSC high-pass filters (\ref{deriv})
for solving Burgers' equation\cite{weijcp99}. 
The DSC kernel parameters 
are chosen as $W=35$ and  ~$\sigma/\Delta=4.5$ in these 
calculations. 
Both $L_{\infty}$ and  $L_{1}$ errors
at 9 different times are listed in TABLE I. 
In an earlier work, 
Kakuda and Tosaka\cite{KakTos} tested their
generalized boundary element method
by using 100 elements, up to 6 iterations and 
the same time increment as ours ($\Delta t=0.01$).
Their results are also listed in TABLE I for a comparison. 
The errors in both methods are very small.       
The DSC results are from ten to $10^5$
times more accurate than those of Kakuda and 
Tosaka\cite{KakTos}(K-T), 
while  they were obtained using much fewer grid points.
Note that at a late time, the accuracy of the present DSC 
algorithm reaches 15 significant figures.

The numerical solution of Burgers' equation at a high 
Reynolds number (Re $= 10^{5}$) is very difficult
due to the presence of shock\cite{JamBon}.
A direct application of the DSC algorithm using 
64 grid points ($N=64$) and a small time increment 
($\Delta t =0.001$) leads
to a  highly oscillatory results as shown in FIG. 2(a).
The time integration is shown up to 0.5 time units and 
eventually collapses at a later time. 
The plot is given in the spatial interval of $[0,2]$, which is  
generated by an antisymmetric extension of the original 
numerical results in the spatial interval of [0,1]. 
The oscillation starts near 0.3 time units
and is accumulated  and amplified in late integrations.
The solid line in FIG. 1 shows the Fourier image of the
result at $t= 0.5$. The image has two large peaks, one near 
the zero and another near the Nyquist frequency ${\pi\over \Delta}$.

To analyze FIG. 2(a) further, a three-scale
wavelet transform is performed
and the result is depicted in  FIG. 2(b). 
Daubechies' biorthogonal 
wavelets (D7/9)\cite{Daubechies} are employed for 
the wavelet transform. 
At the first scale, a response extended 
over a large domain is recorded
over the southeast 
quarter of the quadrangle which corresponds to  
high frequency oscillations in FIG. 2(a). Note that 
high frequency oscillations reside exclusively at the 
southeast quarter of quadrangle because the oscillations 
only occur at a  special frequency range.
The peak in the middle of the high frequency 
response region is due to the shock front, which produces
similar two other narrow high frequency responses at 
the second and third scales respectively. 
The response of the highest amplitude at the 
southwest corner is enlarged in FIG. 2(c). 
Surprisingly, this part seems containing the desired solution 
to Burgers' equation. 
However, there is a kink in the 
late part of the low pass solution 
which does not belong to the desired solution.  
Obviously, had the high frequency oscillations been 
controlled in the course of integration, such a kink 
would not have appeared. 
From this analysis we conclude that the peak near the 
Nyquist frequency in the Fourier image in FIG 1 is due the 
undesired oscillations.

The control of oscillations can be accomplished in a 
number of ways. For example, 
Godunov algorithms\cite{Godunov}, up-wind schemes, 
essentially non-oscillatory (ENO) schemes\cite{HEOC,ShuOsh89}  
and weighted ENO schemes\cite{WENO} 
are standard methods for handling oscillations. 
In the present work, we propose an alternative approach which 
makes use of a conjugated low-pass filter, Eq. (\ref{K1}).  
Since the peak of undesired oscillations is resided 
outside the effective bandwidth of the conjugated
low-pass filter as shown in FIG.1, it can be 
removed by the filtering of the conjugated 
low-pass filter. To eliminate oscillations and preserve
the true solution effectively, we design the following 
conjugated filter oscillation reduction (CFOR) scheme.
We define a wavelet high-pass measure via a multiscale 
wavelet transform of the results of conjugated high-pass 
filter $\{ u(x_k,t_n)\}^N_{k=1}$ (i.e. the numerical solution 
of Burgers' equation) at time $t_n$ as
\begin{equation}\label{wavelet1}
\parallel{\cal W}^{n} \parallel=\sum_m \parallel{\cal W}^n_m \parallel,
\end{equation}
where $\parallel{\cal W}^n_m \parallel$ is given by 
a convolution with a wavelets $\psi_{mj}$ of scale $m$
\begin{equation}\label{wavelet2}
\parallel {\cal W}^n_m \parallel
=\sum_k \left| \sum_j \psi_{mj}(x_k)u^n(x_j) \right|.
\end{equation}   
The CFOR is adaptively implemented 
whenever the high pass measure accesses an appropriate
positive  alarm threshold $\eta$
\begin{equation}\label{wavelet4}
\parallel{\cal W}^{n+1} \parallel-\parallel{\cal W}^{n} \parallel
\geq\eta.
\end{equation}
The choice of $\eta$ depends on the time increment
$\Delta t$ and the grid size $\Delta x$.
FIG. 2(d) shows the results under the same conditions as those
of FIG. 2(a), obtained by using the CFOR scheme.
Note that the oscillations are
eliminated and meanwhile, a sharp shock profile is resolved.
The dots in FIG. 1 shows the Fourier image of the solution at 
time 0.5.  It is noted that there is little change to
the image inside the effective 
bandwidth of conjugated filters. However, peaks of high 
frequency oscillations are effectively eliminated.

We test this scheme for the case of Re$=10^3$ and Re$=10^5$
at $\Delta x=0.01$, $\Delta t=0.001$. 
DSC kernel parameters are the same as stated earlier. 
As shown in FIG. 3(b) and FIG. 3(c), our results are excellent. 
Clearly, all oscillations are effective removed and the 
shock front is very sharp.
In a dramatic case, we consider inviscid Burgers' 
equation (Re=$\infty$). As depicted in FIG. 3(d), 
our CFOR scheme works extremely 
well for this case too. The results for Re$=100$ are shown 
in FIG. 3(a) for a comparison.

In conclusion, a novel approach, the conjugated 
filter oscillation reduction (CFOR) scheme
is introduced for solving Burgers' 
equation with a wide range of Reynolds numbers. 
The essence of the CFOR scheme is to adaptively implement 
a conjugated low-pass filter to effectively remove 
the accumulated numerical errors produced by a set of 
high-pass filters. The conjugated low-pass and high-pass
filters have essentially the same degree of regularity, 
smoothness, time-frequency localization, effective 
support and bandwidth. In this work, all conjugated filters
are constructed by using discrete singular convolution 
kernels\cite{weijcp99}. 

The accuracy of the approach is 
tested at a moderately high Reynolds numbers (Re=100) 
for which analytical solution is available. 
While using much fewer grid points, our results are about 
10 to 10$^{5}$ times more accurate than a previous finite 
element approach\cite{KakTos} .
For extremely high Reynolds numbers (Re$\geq10^5$), 
the DSC algorithm 
develops severe oscillations near the shock
front. The CFOR scheme is 
proposed to minimize error accumulations and to resolve
the shock front. It is found that the present scheme 
is very accurate and robust for integrating Burgers' 
equation over all possible Reynolds numbers.

We note that the present approach is very general.
It can be applied to the numerical solution of other
partial differential equations,
particularly, compressible flows
and hyperbolic conservation laws. 
Moreover, the CFOR scheme can be 
implemented along with any other standard 
computational methods, such as
high-order central difference schemes, finite element methods 
and spectral approximations. A theoretical analysis of the 
present approach and an evaluation of various 
wavelets and filter banks for the CFOR scheme 
will be reported elsewhere\cite{GuWei}.

{This work was supported in part by the National University 
of Singapore. The authors thank Professor C.-W. Shu for 
useful discussions about the ENO and WENO schemes.}

      
        

\newpage
\begin{center}
TABLE I. Comparison of errors for solving  
            Burgers' equation
\end{center}

\begin{center}
\begin{tabular}{c|c|cc} \hline \hline
     & ~K-T~ & \multicolumn{2}{c}{DSC}\\ \cline{2-4}
Time  & $L_{\infty}$ & $L_{\infty}$ & $L_{1}$ \\ \hline
0.4    &~ 2.6(-02) ~&~2.4(-03)  ~&~ 2.2(-04)~\\
0.8    &  2.9(-02)  & 3.3(-03)   &  2.9(-04) \\
1.2    &  1.8(-02)  & 4.7(-04)   &  1.1(-05)\\
3.0    &  6.9(-03)  & 7.6(-08)   &  1.1(-08)\\
10.0   &            & 3.1(-11)   &  1.2(-11)\\
30.0   &            & 1.4(-12)   &  8.6(-13)\\
60.0   &            & 6.9(-14)   &  4.4(-14)\\
90.0   &            & 3.6(-15)   &  2.3(-15)\\
\hline \hline
\end{tabular}
\end{center}

\newpage

\centerline{\bf Figure Captions}

{\bf FIG. 1.} Graph of the frequency responses of the conjugated DSC 
filters (in the unit of $\pi/\Delta$). The maximum 
amplitude is normalized to the unit.
Stars: conjugated low-pass filter; 
Dashed line: 1st order high-pass filter; 
Dash-dots: 2nd order high-pass filter.
Solid line: Fourier image of the numerical solution
of Burgers' equation ($t=0.5$, Re$=10^5$) with oscillations;
Dots: Fourier image of the numerical solution
obtained by using the CFOR scheme.

\vskip 20pt

{\bf FIG. 2.} 
(a) The oscillatory numerical solution of Burgers' equation 
(Re$=10^5, \Delta t=0.001, N=64, t=0\sim 0.5$); 
(b) Three scale wavelet analysis of FIG. 2(a);  
(c) The last scale low pass response of FIG. 2(b); 
(d) The CFOR solution of Burgers' equation (Re$=10^5,
\Delta t=0.001, N=64, t=0\sim 0.5$).

\vskip 20pt

{\bf FIG. 3.} 
The CFOR solutions of  Burger's equation at $t=$0.2 
(i), 0.5 (ii), 1.0 (iii), 1.5 (iv) and 2.0 (v)
($\Delta t=0.001, N=101$).
(a) Re$=100$;
(b) Re$=10^3$;
(c) Re$=10^5$;
(d) Re$=\infty$.

\end{document}